%%% Typos Corrected:  FINAL VERSION 1/25/99
%%% Joint Paper:  Chaatit-Rosenthal
\input amstex
\documentstyle{amsppt}
%\magnification=\magstep1
\NoBlackBoxes
\def\DSC{\operatorname{DSC}}
\def\DBSC{\operatorname{DBSC}}
\def\BDSC{\operatorname{BDSC}}

\def\BPS{\operatorname{BPS}}
\def\PS{\operatorname{PS}}
\def\osc{\operatorname{osc}}
\def\oE{{\overline{E}}}

\def\B{{\Cal B}}
\def\C{{\Cal C}}
\def\U{{\Cal U}}
\def\bU{\text{\bf U}}
\def\real{{\Bbb R}}
\def\nat{{\Bbb N}}
\def\zed{{\Bbb Z}}
\def\th{{\text{th}}}
\def\To{\Rightarrow}
\def\ep{\varepsilon}
\def\chix{{\raise.5ex\hbox{$\chi$}}}
\def\defeq{\mathop{\ \buildrel {\text{\rm def}}\over =\ }\nolimits}
\def\cU{\mathop{\U}\limits^\circ}
\def\Circ#1{\mathop{#1}\limits^\circ}
%%%%%%%%%%%%%%%%%%%%%%%%%%%%%%%%%%%%%%%%%%%%%%%%
\topmatter
\title On Differences of Semi-Continuous Functions\endtitle
\author Fouad Chaatit and Haskell P. Rosenthal\endauthor
\address Fouad Chaatit, School of Science and Engineering, AlAkhawayn
University in Ifrane, Ifrane 53000 Morocco\endaddress
\email F.Chaatit\@alakhawayn.ma\endemail
\address Haskell Rosenthal, Department of Mathematics, The University
of Texas at Austin, Austin, TX 78712-1082\endaddress
\email rosenthl\@math.utexas.edu\endemail
\thanks This research was partially supported by NSF DMS-980153.\endthanks
\abstract Extrinsic and intrinsic characterizations are given for the class
$\DSC(K)$ of differences of semi-continuous functions on a Polish space $K$,
and also decomposition characterizations of $\DSC(K)$ and the class $\PS(K)$
of pointwise stabilizing functions on $K$ are obtained in terms of behavior
restricted to ambiguous sets.
The main, extrinsic characterization is given in terms of behavior
restricted to some subsets of second category in any closed subset of $K$.
The concept of a strong continuity point is introduced, using the
transfinite oscillations $\osc_\alpha f$ of a function $f$ previously
defined by the second named author.
The main intrinsic characterization yields the following $\DSC$ analogue of
Baire's characterization of first Baire class functions:
a function belongs to $\DSC(K)$ iff its restriction to any closed non-empty
set $L$ has a strong continuity point.
The characterizations yield as a corollary that a locally uniformly
converging series $\sum \varphi_j$ of $\DSC$ functions on $K$ converges to
a $\DSC$ function provided $\sum \osc_\alpha \varphi_j$ converges locally
uniformly for all countable ordinals $\alpha$.
\endabstract
\subjclass Primary 26A21,46B03; Secondary 03E15,04A15\endsubjclass
\endtopmatter

\document

\head 1. Introduction\endhead

The interest in functions of the first Baire class can be traced back to
Baire's paper \cite{Ba} in 1899.
In the early twenties S.~Mazurkiewicz and W.~Sierpi\'nski \cite{Mazk},
\cite{Sie} were already studying the subclass of functions that are
the difference of two semi-continuous functions.
{From} a Banach space theory point of view, much of the interest originates
from the $\ell^1$-theorem \cite{R1}.
Subsequently, the goal was to deduce properties of a given separable Banach
space $X$ from the topological class of the functions $f\in X^{**}$.
The fundamental paper \cite{HOR} follows this program.
For example, it was proved in \cite{HOR} that if $K$ is a compact metric
space and $F\in \DSC(K)\sim \C(K)$, then whenever $(f_n)\subset \C(K)$ is a
uniformly bounded sequence converging pointwise to $F$, we have that $c_0$
embeds into $[(f_n)]$, the closed linear span of $(f_n)$.
These efforts culminate in the later result \cite{R2}, characterizing
Banach spaces containing $c_0$.
For further structural results on various subspaces of first Baire class
functions, see also \cite{KL}, \cite{C}, \cite{CMR}, \cite{R3}, \cite{R4}.

We concentrate here on the intrinsic nature of two subclasses of Baire-1
functions: $\PS(K)$, the class of pointwise stabilizing functions, and
$\DSC(K)$, the class of differences of semi-continuous functions on a Polish
space $K$.
The definitions of these classes will be given in the next section.
Note here that the assumption that
$K$ is a Polish space, is a special case of the well-behaved
spaces of W.~Kotze \cite{K}.
It would be interesting to see if the localization theorems presented here
still hold in the more general framework of well-behaved spaces.
We would like to recall here the Baire characterization theorem for a
real-valued function $f$ defined on $K$ to be of the first Baire class in
terms of its restriction to any closed subset $F$ of $K$ having a point of
continuity relative to $F$.
It is in this spirit that we characterize functions of the class $\PS(K)$
in terms of their behavior when restricted to ambiguous sets, i.e., sets
that are simultaneously $F_\sigma$ and $G_\sigma$.
More precisely, we prove in Theorem~3.1 that for a Baire-1 function $F$ on
a Polish space $K$ to be pointwise stabilizing, it suffices that each
$F_{|V^i}$ is continuous, where $(V^i)_{i=1}^\infty$ is a sequence of
ambiguous sets such that $K= \bigcup_{i=1}^\infty V^i$.
Note here that the only if part was already proved in \cite{HOR}.
As corollaries of this, we obtain that a Baire-1 function on a Polish space
with a discrete range is pointwise stabilizing, and moreover such functions
are uniformly dense in the first Baire class.
However, there exists a Baire-1 function on $[0,1]$ whose range is a
convergent sequence, and which is not in $\PS([0,1])$.
In section~4, we prove in Theorem~4.4 that for a bounded real-valued
function $f$ on a Polish space $K$ to be a difference of semi-continuous
function on $K$, it suffices that $f_{|A^j}$ is in $\DSC(A^j)$ for all
$j\ge1$, where $(A^j)_{j=1}^\infty$ is a sequence of ambiguous sets such that
$K=\bigcup_{j=1}^\infty A^j$.
We then use this to give a characterization  of $\DSC$ functions on a Polish
space $K$ in terms of their behavior when restricted to some subset of
second category in closed subsets $L$ of $K$.
More precisely, we prove in Theorem~4.5 that $f$ is in $\DSC(K)$ if and only if
for any closed subset $L$ of $K$ there exists a subset $E$ of $L$ of second
category in $L$ so that $f_{|\oE}$ is in $\DSC(E)$, where $\oE$ is the
closure of $E$.
An advantage of such a characterization is that if $f$ is not in
$\DSC(\Delta)$ where $\Delta$ is the triadic Cantor set, then one can choose
a closed perfect subset $K$ of $\Delta$ so that:
if $f(t) = \sum_{n=1}^\infty f_n(t)$ on $K$, with $f_n$ continuous on $K$
then $\sum_{n=1}^\infty |f_n(t)|< \infty$ only for $t$ in a meager subset
of $K$.
These results deal with what we may call ``localization'':
the theorems say that a function $f:K\to \real$ is in a certain class,
provided its restriction to certain smaller subsets of $K$ are all in
that class.
They are also ``extrinsic''; that is, one needs to know at some level,
that the restrictions of the
function in question are given as  differences of semi-continuous
functions.

The results in Section~5 yield intrinsic characterizations, based on the
invariants introduced in \cite{R2}.
Note that for a given function $f$ on a Polish space $K$, $x\in K$ is a
point of continuity of $f$ iff $\osc f(x) =0$, where $\osc f$ is the
oscillation of $f$.
Letting $\osc_\alpha f$ be the $\alpha^{\th}$ transfinite oscillation
introduced in \cite{R2} (for $\alpha$ a countable ordinal),
we define $x\in K$ to be a {\it strong continuity point of\/} $f$
provided $\osc_\alpha f(x) =0$ for all $\alpha$.
Theorem~5.1 then yields the analogue for $\DSC$ functions of Baire's
famous characterization of Baire class one:
{\it $f$ is $\DSC$ iff $f|L$ has a strong  continuity point for every
non-empty closed subset $L$ of $K$\/}.
The proof of 5.1 also yields an ordinal index for the ``local $\DBSC$
complexity'' of a $\DSC$ function $f$, as well as an effective intrinsic
criterion for determining whether or not a given function is $\DSC$.
Theorem~5.2 yields an answer to the following problem:
{\it if $(f_n)$ is a sequence of $\DSC$ functions converging uniformly to $f$,
what additional properties of the convergence will force $f$ to also be
$\DSC$?}
As noted above, {\it any\/} Baire $-1$ function can be so obtained, without
such conditions.
Theorem~5.2 yields the desired conclusion in terms of the existence of
points of local rapidity of convergence of $(\osc_\alpha f_n|L)$ for
arbitrary non-empty closed subsets $L$ of $K$.
This yields the corollary:
{\it if $\varphi_1,\varphi_2,\ldots$ are $\DSC$ functions with $\sum \varphi_j$
converging locally uniformly to $f$, then $f$ is $\DSC$ provided
$\sum \osc_\alpha \varphi_j$ converges locally uniformly for all
ordinals $\alpha$\/}.

\head 2. Notation\endhead

We use mostly the notation found in \cite{HOR} or \cite{R2}, except that
we do not restrict ourselves to bounded functions; we work only in the field
of real scalars $\real$, although all our results easily extend to the
complex scalars (where $\DSC$ functions are just those whose real and
imaginary parts are $\DSC$.)
Let $K$ be a Polish space.
$\C(K)$ (resp. $C^b(K)$) denotes the class
of continuous (resp. bounded continuous) scalar valued functions on $K$.
$\B_1(K)$ denotes the class of scalar functions on $K$ of the first Baire
class and $\B_1^b (K)$ the class of bounded functions in $\B_1(K)$.
That is,
$$\B_1(K) = \{f:K\to\real :\text{ there exist } (f_n)\text{ in } \C(K)
\text{ with } f_n\to f\text{ pointwise}\}\ .$$
It is easily seen that if $f\in \B_1^b(K)$, then there exists a uniformly
bounded sequence $(f_n)$ in $C^b(K)$ with $f_n\to f$ pointwise.

Recall that $f$ is lower (resp. upper) semi-continuous if for any $\lambda$
the set \break
$\{ x:f(x) >\lambda\}$ is open (resp. $\{x:f(x) \ge \lambda\}$ is
closed).
$f$ is called semi-continuous if it is lower or upper semi-continuous.

This article deals with the following subclass of $\B_1(K)$:
\roster
\item"{}" {\it $\DSC(K) = \{f:K\to\real$: there exists a sequence
$(f_n)\subset \C(K)$ so that $f_n\to f$ pointwise and $\sum_{n=1}^\infty
|f_{n+1}(k) -f_n(k)|<\infty$ for $k\in K$.}
\endroster

We let $\BDSC (K)$ denote the bounded members of $\DSC(K)$.
Again, it is easily seen that such functions may be obtained as the pointwise
limit of uniformly bounded sequences $(f_n)$ in $C^b(K)$ satisfying
$\sum |f_{n+1} (k) -f_n(k)| <\infty$ for all $k$.

By a classical theorem of Baire \cite{H, p.274}, $\DSC(K)$ equals  the
set of all functions that are differences of semi-continuous functions on $K$.

Another subclass is the one of Baire-1 functions on $K$ that are the
difference of two bounded semi-continuous functions on $K$, and will be
denoted by $\DBSC(K)$.
As in \cite{HOR} we will adopt the following equivalent definition:
\roster
\item"{}" $\DBSC(K) := \{f:K\to \real:$ there exist $(\phi_n)$ in $\C(K)$
with $f(k) = \sum_{n=1}^\infty \phi_n(k)$ and $\sum_{n=1}^\infty |\phi_n(k)|
\le C$ for all $k\in K$ and for some constant $C\}$.
\endroster

For $f\in \DBSC(K)$ we then set
$$\align
\|f\|_D  &= \inf \biggl\{ \sup_{t\in K} \sum_{n=0}^\infty |\phi_n(t)|
\ :\ \phi_n\text{ is bounded continuous on  $K$ for all $n$}\\
&\qquad \text{and }\ f(t) = \sum_{n=0}^\infty \phi_n(t),\ t\in K\biggr\}\ .
\endalign$$

As noted in \cite{HOR}, $\DBSC(K)$ is then a Banach algebra under this norm.

\subhead Remarks\endsubhead

1. It is proved in \cite{R2} that for any $\DBSC$ function $f$ on a Polish
space, there exists a countable ordinal $\tau$ so that $f= u-v$, where
$$u = \frac{\|\, |f|+\osc_\tau f\|_\infty - \osc_\tau f+f}2\ ,\
\text{ and }\
v = \frac{\|\, |f|+\osc_\tau f\|_\infty - \osc_\tau f-f}2\ ,$$
$u$ and $v$ being non-negative lower semi-continuous functions; moreover,
the $D$-norm of $f$ is exactly given by $\|f\|_D = \|u+v\|_\infty$.
the transfinite oscillation $\osc_\alpha f$ of $f$ is defined in
Section~5 below.

2. The class $\DBSC(K)$ is in general distinct from the class of Baire-1
functions.
Mazurkiewicz \cite{Mazk} gave a construction of a function that is in
$\B_1(K)\sim \DBSC(K)$ whenever $K$ contains a homeomorphic copy of the
countable compact ordinal $\omega^\omega +1$.
In particular, whenever $K$ is uncountable, the inclusion of $\DBSC(K)$ in
$\B_1(K)$ is strict.
Actually, this result follows via functional analytic reasoning.
Indeed, fix $K$ Polish and assume that $K$ contains a copy of the
(non-compact) ordinal $\omega^\omega$.
Then the uniform and $\DBSC$ norms are {\it inequivalent\/} on $\DBSC(K)$,
hence there exists an $f$ in the {\it uniform closure\/} of $\DBSC(K)$,
not belonging to $\DBSC(K)$, (see \cite{HOR}, \cite{R3}).
Moreover the uniform closure of $\DBSC(K)$, identified as
$\B_{1/2}(K)$ in \cite{HOR}, is also then strictly contained in $\B_1^b(K)$,
and also $\BDSC(K)$ is then strictly larger than $\DBSC(K)$.
Moreover, if $K$ is countable, $\BDSC(K) = \B_1(K) = \ell^\infty (K)$
(see Corollary~3.3 below).

A subclass of $\DSC(K)$ that we also treat here is the class of pointwise
stabilizing functions on $K$, denoted $\PS(K)$.

\proclaim{Definition}
{\it $\PS(K)$ (resp. $\BPS(K)$) is the class of all functions $f:K\to\real$
for which there exists a sequence (resp. uniformly bounded sequence)
of functions $(f_n) \subset C(K)$ so that for any $k\in K$, there is an
integer $m$ with $f_n(k) = f(k)$ for all $n\ge m$.}
\endproclaim
Again, a simple truncation argument yields that $\BPS(K)$ is precisely
the class of bounded functions in $\PS(K)$.

\head 3. A characterization of $\PS(K)$\endhead

\proclaim{Theorem 3.1}
Suppose that $F:K\to\real$ is a given function on a Polish space $K$,
and that $K= \bigcup_{i=1}^\infty W^i$ where each $W^i$ is an ambiguous set.
If $F_{|W^i}$ is continuous for all $i\in \nat$, then $F$ is pointwise
stabilizing.
\endproclaim

\demo{Proof}
By considering $\widetilde{W^j} = W^j \sim \bigcup_{i=1}^{j-1} W^i$
and $\widetilde{W^1} = W^1$ we can suppose without loss of generality
that the $W^i$ are disjoint.
Let
$$ W^j = \bigcup_{n=1}^\infty G_n^j\ ,$$
where for all $n$, $G_n^j$ is closed and $G_n^j\subset G_{n+1}^j$.
Let $f_n: K\to \real$ be a continuous extension of
$F_{|\bigcup_{j=1}^n G_n^j}$.
These $f_n$'s exist by the Tietze extension theorem.
Thus we have.
$$f_{n|G_n^j} = F_{|G_n^j} \text{ for all } 1\le j\le n\ .
\tag 1$$

\proclaim{Claim}
$f_n\to F$ pointwise as $n\to\infty$, and the $f_n$'s stabilize.
\endproclaim

\noindent
Indeed, if $t\in K$ then $t\in W^j$ for some $j$.
Fix such a $j$ and let $n_0\ge j$ large enough so that $t\in G_{n_0}^j$.
Now for all $n\ge n_0$ we have:
$$f_n(t) = f_{n|W^j} (t) = F_{|G_n^j} (t) = F(t)$$
since $t\in G_{n_0}^j$ and (1) holds for all $1\le j\le n_0$.\qed
\enddemo

We note that the result also holds if we assume that the $W^j$ are just
$F_\sigma$'s.
Indeed, for each $i$, choose $A_{ij}$ closed with
$W^i = \bigcup_{j=1}^\infty A_{ij}$.
But evidently $\bigcup_{i,j} A_{ij} = K$ and $f_{|A_{ij}}$ is continuous
for all $i$ and $j$.
Of course, the family $\{A_{ij} :i,j=1,2,\ldots\}$ is countable and each
$A_{ij}$ is ambiguous, being closed: hence Theorem~3.1 applies.
We also note that, since the ambiguous sets form an algebra of sets, the
differences of closed sets are also ambiguous.
We then obtain the following characterization of $\PS (K)$ for Polish
spaces $K$.

\proclaim{Theorem 3.2}
Let $K$ be a Polish space and $F\in \B_1(K)$.
Then the following are equivalent:
\roster
\item"1."  $F\in \PS(K)$.
\item"2." There exists a (disjoint) sequence $(V^i)_{i=1}^\infty$ of
differences of closed subsets of $K$ so that $\bigcup_{i=1}^\infty V^i =K$
and $F_{|V^i}$ is continuous for all $i$.
\item"3." There exists a (disjoint) sequence $(V^i)_{i=1}^\infty$ of
ambiguous subsets of $K$ so that $\bigcup_{i=1}^\infty V^i =K$ and $F_{|V^i}$
is continuous for all $i$.
\endroster
\endproclaim

\demo{Proof}
To see that $1)\To 2)$, recall that from Proposition 4.9 of \cite{HOR},
if $F\in \PS(K)$ then there exist a sequence $(K_n)_{n=1}^\infty$ of closed
sets of $K$ so that $K_n\subset K_{n+1}$, $K= \bigcup_{n=1}^\infty K_n$
and $F_{|K_n}$ is continuous for all $n$.
(The argument in \cite{HOR} does not require that $F$ be bounded.)
It suffices then to set $V^i = K_i\sim K_{i-1}$.\newline
$2)\To 3)$ follows immediately from the comments preceding the statement
of Theorem~3.2.
$3) \To 1)$ follows immediately from the previous theorem.\qed
\enddemo

\proclaim{Corollary 3.3}
Let $K$ be a countable Polish space.
Then every function on $K$ belongs to $\PS (K)$ and hence to $\DSC(K)$.
\endproclaim

\demo{Proof}
Let $F:K\to\real$ be given and let $k_1,k_2,\ldots$ be an enumeration of $K$.
Of course $\{k_i\}$ is ambiguous and $F|\{k_i\}$ is continuous for all $i$,
so $F\in \PS(K)$ by Theorem~3.2 and the obvious fact that
$F\in \B_1 (K)$ (e.g., $F= \sum_n F(k_n) \chix_{k_n}$ pointwise).\qed
\enddemo

\proclaim{Corollary 3.4}
Let $K$ be a Polish space and $F\in \B_1(K)$, so that $F(K)$ is a discrete
set. Then $F\in \PS(K)$.
\endproclaim

\demo{Proof}
The hypothesis means that no point of $F(K)$ is a cluster point of $F(K)$.
It follows easily that $F(K)$ is countable, say $F(K) = \{c_1,c_2,\ldots\}$.
Now a standard result \cite{H} asserts that for $F\in \B_1(K)$, if
$A\subset\real$ is open, then $F^{-1}(A)$ is an $F_\sigma$-set.
For each $i\ge1$, choose $\ep_i>0$ so that $(c_i-\ep_i,c_i+\ep)\cap F(K) =
\{c_i\}$.
Thus $A_i:= F^{-1}\{c_i\} = F^{-1} (c_1 - \ep_i,c_i + \ep_i)$ is an $F_\sigma$.
Since the $A_i$'s are disjoint and partition $K$, it follows that the $A_i$'s
are ambiguous; on the other hand $F_{|A_i}$ is trivially continuous for
all $i$, hence $F\in \PS(K)$ by Theorem~3.1.\qed
\enddemo

\remark{Remark}
In general, no weaker topological assumption on the range of $F$ is possible,
to ensure the validity of Corollary~3.4, when $K$ is uncountable.
Indeed, we may easily, for example, construct a function $f\in \DBSC([0,1])
\sim \PS[0,1]$ whose range is a convergent sequence along with its limit,
as follows:
let $\{d_i\}$ be a countable dense subset of $[0,1]$ and let
$$f= \sum_{j=1}^\infty {\chix_{\{d_j\}} \over 2^j}\ .$$
Then $f$ cannot be pointwise stabilizing, since there is no open non-empty
subset $\U$ of $[0,1]$ so that $f_{|\U}$ is continuous
(cf. Proposition~4.9 of \cite{HOR}).
Evidently $f\in \DBSC([0,1])$, since
$$\|f\|_{\DBSC} \le \sum_{j=1}^\infty
{\|\chix_{ d_j}\|_{\DBSC}\over 2^j} \le 2\ .$$
Finally, $f([0,1]) = \{0\} \cup \{2^{-j} :=1,2,\ldots\}$, a
convergent sequence with its limit.
\endremark

Following Hausdorff \cite{H}, let us say that $F\in \B_1(K)$ is a
{\it step function\/} if $F(K)$ is discrete.

\proclaim{Corollary 3.5}
Let $K$ be a Polish space and $F\in \B_1(K)$.
Then $F$ is a uniform limit of step functions.
Hence $\DSC (K)$ is uniformly dense in $\B_1(K)$.
\endproclaim

\demo{Proof}
We first recall the classical

\example{Fact 1}
Given $G\subset F\subset K$, $G$ a $G_\delta$ and $F$ an $F_\sigma$, there
exists an ambiguous set $A$ with $G\subset A\subset F$
(cf. \cite{Ku}, Theorem 2, page 350).
\endexample

\example{Fact 2}
$f:K\to\real \in \B_1(K) \Leftrightarrow f^{-1}(U)$ is an $F_\sigma$
for all open $U\subset \real$ (see \cite{H}).
\endexample

Now fix $n\ge 2$ and for each $m\in\zed$, $1\le j\le n$, choose
$A_j^m$ an ambiguous set with
$$G_j^m \defeq f^{-1} \left(\Big[ \frac{m+j-1}n, \frac{m+j}n \Big]\right)
\subset A_j^m \subset f^{-1}
\left( \frac{m+j-2}n, \frac{m+j+1}n\right) \defeq F_j^m\ .
\tag 2$$
(This is possible by the Facts, since $G_j^m$ is then a $G_\delta$,
$F_j^m$ an $F_\delta$.)
Then evidently
$$\bigcup_{m=-\infty}^\infty \bigcup_{j=1}^n A_j^m = K\ .
\tag 3$$
Finally, by disjointifying the $A_j^m$'s, choose a sequence $W_1,W_2,\ldots$
of disjoint ambiguous sets with $K= \bigcup_{j=1}^\infty W_j$ so that for
all $i$, there is an $m(i) \in\zed$ and $j(i)$, $1\le j(i)\le n$ with
$$f(W_i) \subset \left( \frac{m(i) + j(i)-2}n, \frac{m(i)+j(i)+1}n\right)\ .
\tag 4$$
Then let
$$f_n = \sum_{i=1}^\infty \left( \frac{m(i)+j(i)}n\right)
\chix_{W_i}\ .
\tag 5$$
Now we have by Theorem 3.1 that $f_n\in \PS (K)$; of course $f_n$ is thus
a step function.
Now for any $i$ and $x\in W_i$, $|f(x)-f_n(x)| \le \frac3n$.
But then $\|f- f_n\|_\infty \le \frac3n$, whence $f_n\to f$ uniformly.\qed
\enddemo

\remark{Remarks}
1. Corollaries 3.3 and 3.5 again reveal the immense difference between
$\DBSC(K)$ and $\BDSC(K)$.
Indeed, for general $K$, the uniform closure of $\DBSC(K)$ equals
$\B_{1/2} (K)$, a very thin subset of $\B_1^b (K)$ as long as $K$ contains
a homeomorphic copy of $\omega^\omega$.
Of course Corollary~3.5 immediately yields for any Polish $K$, that
the uniform closure of $\BDSC(K)$ equals $\B_1^b(K)$.
Moreover for $K$ as above, $\DBSC(K)\subsetneqq \B_{1/2}(K)$, but if $K$ is
countable,
$$\BPS(K) = \BDSC(K) = \B_1^b (K) = \ell^\infty (K)$$
by Corollary 3.3.

2. Corollaries 3.3--3.5 are essentially given (with different reasoning)
in \cite{H} (see page 278).
\endremark

\head 4. Extrinsic characterizations of $\DSC(K)$\endhead

To motivate our theorems, we first establish the following:

\proclaim{Proposition 4.1}
Let $K$ be a perfect Polish space.
If $U$ is a non-empty open subset of $K$ and $\lambda\in\real$ then the
set $\B_{\U,\lambda} := \{f\in \B_1^b(K):$ there exist $(\phi_n)\subset
C^b(K)$ with $\sum_{n=1}^\infty \phi_n =f$ and $\sum_{n=1}^\infty |\phi_n(t)|
\le\lambda$ for all $t\in \U\}$ is of first category in $\B_1^b (K)$.
\endproclaim

\demo{Proof}
Let $X = \{f\in \B_1^b (K) : f_{|\U} \in \DBSC(\U)\}$.
Then $X$ is a linear subspace of $\B_1(K)$; we may introduce a norm on $X$ by
$$\|x\| = \|x_{|\U}\|_{\DBSC(\U)} + \|x_{|\U^c}\|_\infty\ .$$
Now if $Y = \{f\in \DBSC (K) : f\equiv 0 $ on $\U^c\}$, then $Y$ is a linear
subspace of $X$, and the norms $\|\cdot\|$ and $\|\cdot\|_\infty$ are not
equivalent on $Y$.
Indeed, since $\U$ is open, $\U$ is again a perfect Polish space, and hence
uncountable; it is trivial that $Y$ is canonically isometric to
$\DBSC(\U)$, so this follows from \cite{HOR}.
Hence, the norms $\|\cdot\|$ and $\|\cdot\|_\infty$ are not equivalent on $X$.
Now evidently $\B_{\U,\lambda}\subset X$, so the proposition follows
from the easy:
\enddemo

\proclaim{Lemma 4.2}
If $X$ is a linear subspace of a Banach space $B$ so that $\|\cdot\|_B
\le \|\cdot\|_X$ where $\|\cdot\|_X$ is a norm on $X$ which is not equivalent
to the norm $\|\cdot\|_B$ of $B$, then $X$ is of first category.
\endproclaim

\demo{Proof}
Let $F_n: = \overline{\{x\in X:\|x\|_X \le n\}}$ (where the closure is
taken in $B$).
Then $F_n$ has void interior, for if not, by the standard proof of the
open mapping theorem, the norms are equivalent on $X$ (and $X=B$!).
So $F_n$ is meager.
Hence $\bigcup_{n=1}^\infty F_n$ is of first category.
But $X\subset \bigcup_{n=1}^\infty F_n$, so $X$ is of first category.\qed
\enddemo

\proclaim{Corollary 4.3}
If $K$ is a perfect Polish space, then $\BDSC (K)$ is of first
category in $\B_1^b(K)$.
\endproclaim

\demo{Proof}
Let $f\in \BDSC(K)$ and $(\phi_n)_{n=1}^\infty \subset C^b(K)$ with
$f= \sum_{n=1}^\infty \phi_n$ and $\sum_{n=1}^\infty |\phi_n (t)| < \infty$
for all $t\in K$.
Let $K_m := \{t\in K: \sum_{n=1}^\infty |\phi_n (t)| \le m\}$.
Then $K = \bigcup_{m=1}^\infty K_m$.
But $K$ is a Baire space, so there exist $m_0 \in\nat$ and $\U$ a non-empty
open subset of $K$ with $\U\subset K_{m_0}$.
So $f\in \B_{\U,m_0}$ and therefore, if $(\U_n)_{n=1}^\infty$ is a basis of
neighborhoods for $K$, then:
$$\BDSC(K) \subset \bigcup_{n,m} \B_{\U_n,m} \ ,$$
which proves that $\BDSC(K)$ is of first category.\qed
\enddemo

\remark{Remarks}
1. It also follows that if $K$ is a perfect Polish space, then $\DSC(K)$
is of first category in $\B_1(K)$, endowed with the topology of
uniform convergence.
($\B_1(K)$ is a complete metric space in this topology, where  e.g.,
we set
$\rho (f,g) = \sup_{x\in K} \frac{|f(x)-g(x)|}{1+ |f(x)-g(x)|}$.)
To see this, let $\U$ be a non-empty open subset of $\B_1(K)$ and let
$D_{\U} = \{f\in \B_1(K): f\mid \U\in \DBSC(\U)\}$.
Now since $\chix_{\U}$  itself is in $\B_1^b(K)$, we have that
$\B_1(K) = \B_1(\U) \oplus \B_1(\sim \U)$ via the obvious identifications.
But then $D_{\U} = \DBSC(\U) \oplus \B_1(\sim \U)$, and so $\DBSC(\U)$
is first category in $\B_1^b(\U)$ by Corollary~4.3.
But $\B_1^b(\U)$ is closed in $\B_1(\U)$, so $\DBSC(\U)$ is first
category in $\B_1(\U)$, whence $D_{\U}$ is first category  in $\B_1(K)$.
But now the argument for 4.3 goes through (deleting the
``$\vphantom{\B}^b$'' superscripts) yielding that
$\DSC(K)\subset \bigcup_{n=1}^\infty D_{\U_n}$, hence $\DSC(K)$ is
first category.

2. For any Polish space $K$, let $\B_{1/2}(K)$ denote the {\it uniform
closure\/} of $\DBSC(K)$ in $\B_1^b (K)$.
An intrinsic equivalent definition may be found in \cite{HOR}, where it is
shown that $\B_{1/2}(K)\ne \B_1^b (K)$ and $\B_{1/2} (K)\ne \DBSC(K)$ if
$K$ contains a subset homeomorphic to $\omega^\omega +1$
(see Proposition~5.3 of \cite{HOR}).
It then follows by the same argument as above that
{\it if $K$ is a perfect Polish space, then $\B^{1/2}(K) \cap \BDSC (K)$
is of first category\/}.
For the proof, simply replace ``$\B_1^b(K)$'' by ``$\B_{1/2}(K)$'' in
Proposition~4.1 and its proof, and ``$\BDSC(K)$'' by
``$\BDSC(K)\cap \B_{1/2}(K)$'' in the proof of 4.3.
\endremark

The proof of Corollary 4.3 actually shows that $\DSC(K)$ is a subset of the
class of functions $f:K\to\real$ so that, for any closed subset $L$ of $K$,
there exists a relatively open subset $\U$ of $L$ such that $f_{|\U}$ is
in $\DBSC(\U)$  (since for any closed subset $L$ of $K$ we obviously have
that $\DSC (K)_{|L}\subset \DSC(K)$).
We will prove in our second intrinsic characterization theorem that this
inclusion is in fact an equality.
This theorem requires the following crucial decomposition result (our first
intrinsic characterization).

\proclaim{Theorem 4.4}
Let $K$ be a Polish space and $f: K\to\real$ be a given function.
Suppose that $K= \bigcup_{j=1}^\infty A^j$, where $A^j$ is ambiguous and
$f_{|A^j} \in \DSC(A^j)$ for all $j\ge1$.
Then $f\in \DSC(K)$.
\endproclaim

\demo{Proof}
Without loss of generality, we can suppose that the $A^j$'s are disjoint.
By hypothesis,
$$f_{|A^j} = \lim_{n\to\infty} f_n^j \text{ with }
\sum_{n=1}^\infty |f_{n+1}^j (t) - f_n^j (t)| < \infty
\text{ for all } t\in A^j\ ,
\tag 6$$
where $f_n^j$ is continuous on $A^j$.
Now each $A^j$  is an $F_\sigma$, say $A^j = \bigcup_{n=1}^\infty A_n^j$
with the $A_n^j$'s closed and in addition, we can suppose that the $A_n^j$'s
are increasing in $n$, i.e., $A_n^j \subset A_{n+1}^j \subset \cdots$,
while of course they are disjoint in $j$.
The Tietze extension theorem then provides for each $n$ a continuous
function $f_n :K\to\real$ that extends $\sum_{j=1}^n f_n^j \chix_{A_n^j}$.
In other words we have:
$$f_{n|A_n^j} = f_{n|A_n^j}^j\ \text{ for all }\ 1\le j\le n\ .
\tag 7$$

\remark{Claim}
$f_n\to f$ pointwise and $\sum_{n=1}^\infty |f_{n+1} (t) - f_n(t)| <\infty$
for all $t$.
\endremark
\medskip

\noindent
Indeed, let $t\in K$; then $t$ lies in some $A^j$.
Pick $n_0\ge j$ large enough so that $t\in A_{n_0}^j$.
Now, if $n\ge n_0$ then $f_n(t) = f_n^j (t)$ by (7) since $n_0 \ge j$,
and $f_n^j (t)\to f(t)$ as $n\to\infty$ by (6).
Also for that same $n_0$,
$$\sum_{n\ge n_0} |f_{n+1} (t) - f_n(t)|
=  \sum_{n\ge n_0} |f_{n+1}^j (t) - f_n^j (t)|
< \infty$$
by (6).\qed
\enddemo

We are now ready for our second intrinsic characterization theorem.

\proclaim{Theorem 4.5}
Let $f:K\to\real$ be a given bounded function on some Polish space $K$.
Then the following are equivalent:
\roster
\item"1." $f\in \DSC(K)$.
\item"2." For any closed non-empty subset $L$ of $K$, there exists a closed
relative neighborhood $\U$ in $L$ so that $f_{|\U}$ is in $\DSC(\U)$.
\item"3." For any closed non-empty subset $L$ of $K$, there exists a closed
relative neighborhood $\U$ in $L$ so that $f_{|\U}$ is in $\DBSC(\U)$.
\item"4." For any closed non-empty subset $L$ of $K$, there exists a subset
$E$ of second category in $L$ so that $f_{|\oE}$ is in $\DSC(E)$.
\endroster
\endproclaim

Recall that $\U\subset L$ is a relative neighborhood in $L$ if $\U$ has
non-empty relative interior with respect to $L$.

\remark{Remark}
Condition 4. means that for any closed subset $L$ of $K$ there exist a
subset $E$ of second category in $L$ and a sequence of continuous functions
$(f_j)$ on $\oE$, so that $(f_j)$ converges to $f$ pointwise on $\oE$
and $\sum_{j=1}^\infty |f_{j+1} (t) - f_j(t)| <\infty$ for all $t\in E$.
\endremark

\demo{Proof}
We will prove that:  $(1)\To (3)\To (2) \To (1)\To (4) \To (3)$.

$(1)\To (3)$:
Let $(\phi_n)_{n=1}^\infty \subset \C(K)$ with $f= \sum_{n=1}^\infty \phi_n$
and $\sum_{n=1}^\infty |\phi_n (t)|<\infty$ for all $t\in K$.
Let $L$ be a closed subset of $K$.
Then $L = \bigcup_{m=1}^\infty L_m$ where $L_m:= \{t\in L:\sum_{n=1}^\infty
|\phi_n(t)| \le m\}$.
Then $L_m$ is closed for all $m$, and so there exists, by the Baire category
theorem, an integer $m_0$ so that the interior of $L_{m_0}$ is not empty.
Take $\U$ to be the closure of that interior.
Then $\U$ is a closed subset of $L$ with a non-empty interior and
$f_{|\U}$ is in $\DBSC (\U)$.

$(3)\To (2)$ is trivial.

To prove $(2)\To (1)$ we start by taking $L_0 =K$.
By hypothesis there exists $\U_0$ a closed relative neighborhood in $L_0$
so that $f_{|\U_0}$ is in $\DSC (\U)$.
Let $v_0 := \cU{\!}_0$ and set $L_1 := L_0\sim V_0$.
Then $L_1$ is closed and consequently there exists a closed relative
neighborhood $\U_1$ so that $f_{|\U_1}$ is in $\DSC(\U_1)$.
Let $V_1 := \cU{\!}_1$ and set $L_2 := L_1\sim V_1$.
Again $L_2$ is closed and so on.
We thus construct a decreasing family of closed sets $(L_\alpha)_\alpha$
with $L_{\alpha+1} := L_\alpha \sim V_\alpha$; and at limit ordinals,
$L_\alpha : = \bigcap_{\beta<\alpha} L_\beta$.
So, $V_\alpha = L_\alpha \sim L_{\alpha +1}$, and
$$f_{|V_\alpha } = f_{|L_\alpha \sim L_{\alpha+1}} \in \DSC (V_\alpha)\ .
\tag 8$$
Since $K$ is metrizeable separable, $L_\eta = \emptyset$ for some
$\eta <\omega_1$.
Now
$$K = \bigcup_{\alpha<\eta} L_\alpha \sim L_{\alpha+1}\ .
\tag 9$$
To see this, take a $t\in K= L_0$.
If $t$ is not in $L_1$, we are done; otherwise, let $\beta:=\sup\{\alpha :
t\in L_\alpha\}$.
Then $t\in L_\beta \sim L_{\beta+1}$.
Now, since $\eta$ is countable, we can enumerate the ordinals that are less
than $\eta$, as $\beta_1,\beta_2,\ldots$.
Thus, setting $M_n= L_{\beta_n} \sim L_{\beta_n+1}$ for all $n$, we have that
$$K= \bigcup_{n=1}^\infty M_n\ ,
\tag 10$$
and each set $M_n$ is ambiguous, begin a difference of closed sets.
Combining Theorem~4.4 with equations (8) and (10), we get that $f$ must be
in $\DSC(K)$.

$(1)\To (4)$ is trivial.

$(4)\To (3)$:
Let  $L$ be a closed set in $K$.
By hypothesis, there exists a subset $E$ of second category in $L$ so that
$f_{|\oE}$ is in $\DSC(E)$.
Let then $(f_n)\subset \C(\,\oE)\,$ with $f_{|\oE} (t) = \lim_{n\to\infty}
f_n(t)$ for all $t\in \oE$, and $\sum_{n=1}^\infty |f_{n+1} (t)-f_n(t)|
<\infty$ for all $t\in E$.
Now let $F: = \{t\in \oE: \sum_{n=1}^\infty |f_{n+1}(t) -f_n(t)|<\infty\}$.
Then obviously $E\subset F$.
But $F= \bigcup_{m=1}^\infty F_m$ where $F_m : = \{t\in \oE: \sum_{n=1}^\infty
|f_{n+1}(t) -f_n(t)| \le m\}$.
$F_m$ is easily seen to be closed; and, since $E$ is of second category
in $L$, there exists an integer $m_0$ so that $\Circ{F}{\!}_{m_0} \ne
\emptyset$.
Setting then $\U = F_{m_0}$ ends the proof.\qed
\enddemo

\remark{Remark}
The following is another equivalent statement in Theorem 4.5:
\roster
\item"$4'$." {\it For any closed subset $L$ of $K$ there exists an $F_\sigma$
subset $E$ of $L$ of second category in $L$ so that $f_{|E}$ is in
$\DSC(E)$.}
\endroster
\endremark

\noindent
To see this simply  let $E= \bigcup_{m=1}^\infty E_m$ where $E_m$ are closed
subsets of $L$.
Since $E$ is of second category in $L$, some $E_m$ has a non-empty interior
relative to $L$.
But of course, if $f_{|E}$ is in $\DSC(E)$, then
$f_{|E_m} \in \DSC (E_m)$.

\head 5. Intrinsic characterizations of $\DSC(K)$\endhead

We first recall the transfinite oscillations $\osc_\alpha f$ of a given
function $f$ defined on a Polish space $K$, as introduced in \cite{R2}.
For any extended real valued function $g$ on $K$ and $x\in K$,
$\varlimsup_{y\to x} g(y)$ denotes the ``unrestricted'' lim~sup of $g$ as
$y$ tends to $x$:
$\varlimsup_{y\to x} g(y) = \inf_U \sup g(U)$, the inf over all open
neighborhoods $U$ of $x$.
$\bU g$ denotes the upper semi-continuous envelope of $g$:
for $x\in K$, $\bU g(x) = \varlimsup_{y\to x} g(y)$.

\proclaim{Definition}
The $\alpha^{\th}$ oscillation of $f$, $\osc_\alpha (f)$, is defined by
ordinal induction as follows:
set $\osc_0 f\equiv 0$.
Suppose $\beta >0$ is a given ordinal and $\osc_\alpha f$ has been defined
for all $\alpha <\beta$.
If $\beta$ is a successor, say $\beta = \alpha+1$, we define
$$\widetilde{\osc}_\beta f(x) = \varlimsup_{y\to x} (|f(y) - f(x)|
+ \osc_\alpha f(y))\ \text{ for all }\ x\in K\ .
\tag 11$$
If $\beta$ is a limit ordinal, we set
$$\widetilde{\osc}_\beta f = \sup_{\alpha <\beta} \osc_\alpha f\ .
\tag 12$$
Finally, we set $\osc_\beta f = \bU \widetilde{\osc}_\beta f$.
\endproclaim

Evidently $\osc_\alpha f$ is a $[0,\infty]$-valued upper semi-continuous
function for all $\alpha$.
A motivating result for the following:
Theorem~3.5 of \cite{R2} yields that $f$
{\it is locally in $\DBSC(K)$ iff $\osc_\alpha f$ is real-valued for
all $\alpha<\omega_1$}.

Classically, $\osc f$ is defined by the equation
$$\osc f(x) = \varlimsup_{y,z\to x} |f(y)- f(z)|
= \varlimsup_{y\to x} f(y) - \varliminf_{y\to x} f(y)
\tag 13$$
(where in the last identity, the meaningless term ``$\infty-\infty$''
is replaced by ``$\infty$'' if it occurs).

Now we easily have that $\osc_1 f \le \osc f\le 2\osc_1f$; thus $f$
is continuous at $x$ iff $\osc_1 f(x)=0$.

\proclaim{Definition}
$x\in K$ is called a strong continuity point of $f$ if
$\osc_\alpha f(x) = 0$ for all $\alpha$.
\endproclaim

Our first intrinsic characterization theorem yields a $\DSC$ analogue of
Baire's theorem characterizing  Baire-1 functions.

\proclaim{Theorem 5.1}
Let $K$ be a Polish space and $f:K\to \real$ a given function.
Then the following are equivalent:
\roster
\item"1." $f\in \DSC(K)$.
\item"2." For any closed non-empty subset $L$ of $K$, $f|L$ has a strong
continuity point.
\item"3." For any closed non-empty subset $L$ of $K$, the set of strong
continuity points of $f|L$ contains a dense $G_\delta$ subset of $L$.
\item"4." For any closed non-empty subset $L$ of $K$, there exists an $x\in L$
so that $\sup_\alpha \osc_\alpha f(x) <\infty$.
\endroster
\endproclaim

\demo{Proof}
We show $1\To 3\To 2\To 4\To 1$.
Of course $3\To 2\To 4$ are trivial.

$1\To 3$:
Let $f= u-v$ with $u$ and $v$ upper semi-continuous.
Now by Baire's famous theorem, if $U,V$ denote the set of points of
continuity of $u|L$ and $v|L$ respectively, $U$ and $V$ are both dense
$G_\delta$'s, hence $G \defeq U\cap V$ is a dense $G_\delta$ subset of $L$.
Now if $x\in G$, then $\osc u|L (x) = \osc v|L(x)=0$.
But for any $\alpha <\omega_1$, $\osc_\alpha u|_L = \osc u|_L$ and
$\osc_\alpha v|_L = \osc v|_L$ since $u,v$ are semi-continuous, by
Proposition~3.4(d) of \cite{R2}.
Hence we have that
$$\align
\osc_\alpha f|L (x) &= \osc_\alpha (u-v)|L(x) \le \osc_\alpha u|L (x)
+ \osc_\alpha v|L (x) \tag 14\\
&\hskip1in \text{ (by Proposition 3.4(b) of \cite{R2})}\\
&= 0\ .\endalign $$
Thus $\osc_\alpha f|L(x) =0$ for all $\alpha$, so $x$ is a strong point of
continuity of $f|L$.

$4\To 1$:
Let $L$ be a closed non-empty subset of $K$.
By Lemma~3.7 of \cite{R2} there exists an $\eta <\omega_1$ so that
$\osc_\eta f|_L = \osc_\beta f|_L$ for all $\beta >\eta$.
Let $U = \{x\in L: \osc_\eta f|L <\infty\}$.
Then $U$ is non-empty since we assume 4 holds.
But $U$ is a relatively open subset of $L$ since $\osc_\eta f|L$ is
upper semi-continuous.
But then by Proposition~3.4(c) of \cite{R2},
$(\osc_\alpha f|_L \pm f|_L)$ are both upper semi-continuous, hence
$f|L\in \DSC(U)$.
Thus condition~3 of Theorem~4.5 holds, whence $f\in \DSC(K)$ by its
conclusion.\qed
\enddemo

\remark{Remark}
The proof of this result yields an ``effective'' ordinal index for $\DSC$
functions on a Polish space $K$ (as well as an ``effective'' intrinsic
criterion for determining {\it when\/} a function is $\DSC$), as follows.
Following \cite{R2}, we define, for a general $f:K\to\real$, $i_D(f)$,
the $D$-index of $f$, by $i_D(f) = \min \{\alpha<\omega_1: \osc_\alpha f
= \osc_{\alpha+1}f\}$.
Now let $\eta_1 = \eta_1 (f) = i_D(f)$ and
$K_1 (f) = \{ x:\osc \eta_1 (f) =\infty\}$.
(Also set $K_0=K$.)
It {\it follows\/} from our argument for 5.1 (i.e., the cited results in
\cite{R2}) that $K\sim K_1 = \{x\in K:f$ is locally $\DBSC$ at $x\}$
(where $f$ is locally $\DBSC$ at $x$ if  $\exists\ U$ a neighborhood of $x$
with $f|U \in \DBSC(U)$).
Of course $K_1$ is closed (possibly empty).
Now for each ordinal $\beta<\omega_1$ having defined $K_\beta$ for all
$\alpha <\beta$ set $K_\beta = \bigcap_{\alpha<\beta} K_\alpha$ if $\beta$
is a limit ordinal.
If not, let $\alpha +1 =\beta$; if $K_\alpha =\emptyset$, set
$K_\beta =\emptyset$.
Otherwise, let $\eta_\beta = \eta_\beta (f) = i_D (f|K_\alpha)$ and let
$K_{\alpha+1} = \{x\in K_\alpha : \osc_{\eta_\beta} f|K_\alpha (x)
= \infty\}$.
Now since $K$ is Polish, we may {\it define\/} $i_{\DSC} (f) $ to be
the least $\alpha <\omega_1$ so that $K_\alpha = K_{\alpha+1}$.
{\it It then follows that $f\in \DSC (K)$ iff $K_\alpha = \emptyset$,
where $\alpha = i_{\DSC}(f)$.}
Indeed, when $K_\alpha =\emptyset$, $K=\bigcup_{\beta<\alpha} K_\beta \sim
K_{\beta+1}$ and $f|(K_\beta\sim K_{\beta+1}) \in \DSC(K_\beta\sim
K_{\beta+1})$ for all $\beta$, whence $f\in \DSC(K)$,
by Theorem~4.4.
In this case, the index $i_{\DSC}(f)$ measures the
``$\DBSC$-complexity of the function $f$.
It seems very likely that if $K$ is uncountable, then there exist
$\BDSC$ functions on $K$ with arbitrarily large countable $\DSC$ indexes.
\endremark

Our last result yields an answer to the following problem:
Suppose $f_n\to f$ uniformly on $K$, with $f_n$ in $\DSC(K)$ for all $n$.
What additional properties of the convergence of the sequence $(f_n)$
will force $f$ to also belong to $\DSC(K)$?

\proclaim{Theorem 5.2}
Let $K$ be a Polish space, and let $f,f_1,f_2,\ldots$ be real-valued functions
on $K$.
Suppose for all $\alpha <\omega_1$ and non-empty closed subsets $L$ of $K$,
there exists a subsequence $(f'_j)$ of $(f_j)$ and a closed relative
neighborhood $U$ of $L$ so that
\roster
\item"(i)" $f'_j\to f$ uniformly on $U$
\endroster
and
\roster
\item"(ii)" $\sum_j \|\osc_\alpha (f'_{j+1} - f'_j)|U)\|_\infty <\infty$.
\endroster
Then $f \in \DSC (K)$.
\endproclaim

We first require the following result from \cite{R4}, which for the sake
of completeness we prove again here.

\proclaim{Lemma 5.3}
Let $\alpha <\omega_1$,  $K$ a Polish space, and $\varphi,\varphi_1,
\varphi_2,\ldots$ real-valued functions on $K$ so that
\roster
\item"(i)" $\sum_{j=1}^n \varphi_j\to\varphi$ uniformly on $K$
\endroster
and
\roster
\item"(ii)" $\sum \|\osc_\alpha \varphi_j\|_\infty <\infty$.
\endroster
Then $\|\osc_\alpha (\varphi - \sum_{j=1}^n \varphi_j)\|_\infty \to0$
as $n\to\infty$.
\endproclaim

\remark{Remark}
It follows immediately that $\|\osc_\alpha\varphi\|_\infty <\infty$.
Indeed, for all $n$,
$$\align \|\osc_\alpha \varphi\|_\infty
&\le  \Big\|\osc_\alpha
\biggl( \varphi -\sum_{j=1}^n \varphi_j\biggr)\Big\|_\infty
+ \Big\| \osc_\alpha \sum_{j=1}^n \varphi_j\Big\|_\infty \\
&\le \big\| \osc_\alpha \biggl( \varphi -\sum_{j=1}^n \varphi_j\biggr)
\Big\|_\infty
+ \sum_{j=1}^\infty \|\osc_\alpha \varphi_j\|_\infty\ .
\endalign$$
Hence in fact
$\|\osc_\alpha \varphi\|_\infty
\le \sum_{j=1}^\infty \|\osc_\alpha \varphi_j\|_\infty <\infty$.
\endremark

\demo{Proof of Lemma 5.3}
Let $f_n = \sum_{j=1}^n \varphi_j$ for all $n$.
Also let $f_0=0$.
We prove by induction on $\gamma \le \alpha$ that
$$\osc_\gamma (\varphi -f_n) (x)
\le \sum_{j=n+1}^\infty (\osc_\gamma \varphi_j) (x)\ \text{ for all }
x\in K\ ,\ \text{ all } n= 0,1,2,\ldots\ .
\tag 15$$

Of course (15) yields in particular that
$$\|\osc_\gamma (\varphi -f_n)\|_\infty
\le \sum_{j=n+1}^\infty \|\osc_\gamma \varphi_j\|_\infty \to 0
\ \text{ as }\ n\to\infty\ ,$$
since $\osc_\gamma g\le \osc_\alpha g$ for any function $g$.

Suppose then $0\le\gamma <\alpha$ and (15) has been proved for $\gamma$.
Now fixing $n$ and $\ep >0$, choose  $q>n$ so that
$$\gather
\Big| (\varphi -f_n)(y) - \sum_{j=n+1}^q \varphi_j(y)\Big| \le\ep
\ \text{ for all } \ y\in K
\tag 16\\
\intertext{and}
\sum_{j=q+1}^\infty \|\osc_\alpha \varphi_j\|_\infty \le \ep \ .
\tag 17
\endgather$$
Now, fixing $x\in K$, we have for any $y\in K$ that
$$\align
&|(\varphi - f_n)(y) - (\varphi - f_n)(x)| + \osc_\gamma (\varphi -f_n)(y)\\
&\qquad \le \Big| \sum_{j=n+1}^q \varphi_j (y)-\varphi_j(x)\Big|
+ \sum_{j=n+1}^\infty \osc_\gamma \varphi_j (y) + 2\ep \\
&\hskip1.5truein \text{ (by (16) and the induction hypothesis (15)}\\
&\qquad \le \sum_{j=n+1}^q \left( |\varphi_j(y) -\varphi_j(x)|
+ \osc_\gamma \varphi_j(y)\right) +3\ep \\
&\hskip1.5truein \text{ by the triangle inequality and (17).}
\endalign$$
Thus by definition,
$$\align
\widetilde{\osc}_{\gamma+1} (\varphi -f_n)(x)
& \le \varlimsup_{y\to x} \sum_{j=n+1}^q \left( |\varphi_j(y) -\varphi_j(x)|
+ \osc_\gamma \varphi_j(y)\right)
\tag 18\\
& \le \sum_{j=n+1}^q \widetilde{\osc}_{\gamma+1} \varphi_j(x) + 3\ep\\
& \le \sum_{j=n+1}^q \osc_{\gamma+1} \varphi_j (x) + 3\ep\ .
\endalign$$
Now since $\sum_{j=n+1}^q \osc_{\gamma+1} \varphi_j$ is upper semi-continuous
and $x$ is an arbitrary point in $K$,
$$\osc_{\gamma+1} (\varphi -f_n) \le \sum_{j=n+1}^\infty
\osc_{\gamma+1} \varphi_j + 3\ep \ \text{ pointwise.}
\tag 19$$
Of course since $\ep>0$ is arbitrary, (15) is established for $\gamma+1$.

Finally, suppose $\beta\le\alpha$ is a limit ordinal and (15) is established
for all $\gamma <\beta$.
But then fixing $x\in K$, we have
$$\align
\widetilde{\osc}_\beta (\varphi -f_n)(x)
= \sup_{\gamma<\beta} \osc_\gamma (\varphi -f_n)(x)
& \le \sup_{\gamma <\beta} \sum_{j=n+1}^\infty \osc_\gamma \varphi_j(x)
\tag 20\\
& \le \sum_{j=n+1}^\infty \osc_\beta \varphi_j (x)\ .
\endalign$$
Now again by taking the upper semi-continuous envelope, we obtain
from (20) that (15) holds for $\gamma =\beta$,
completing the proof of the lemma by transfinite induction.\qed
\enddemo

\remark{Remark}
Of course the lemma yields that if $X_\alpha (K) = X_\alpha$  is the class
of bounded functions $f$ on $K$ with $\osc_\alpha f$ bounded, then $X_\alpha$
is a Banach space under the norm $\|f\|_{X_\alpha} =  \max \{\|f\|_\infty,
\|\osc_\alpha f\|_\infty\}$.
In fact the $X_\alpha (K)$'s are Banach algebras, with a rich structure
connected with invariants for general Banach spaces; see \cite{R4}.
\endremark

We are finally prepared for the

\demo{Proof of Theorem 5.2}
Let $L$ be a non-empty closed subset of $K$, and (as in the proof of
Theorem~5.1), choose $\alpha$ a countable ordinal so that $\osc_\alpha f|_L
= \osc_\beta f|_L$ for all $\beta > \alpha$.
Now choose $U$ and $(f'_j)$ as in the hypotheses of Theorem~5.2.
Then setting $\varphi_j = f'_{j+1} - f'_j$ for all $j>1$, $\varphi_1=f'_1$,
we obtain by Lemma~5.3 that $\|\osc_\alpha (f-f'_n)|U\|_\infty \to0$
as $n\to\infty$, whence (by the remark following 5.2),
$\|\osc_\alpha f|U\|_\infty <\infty$.
It then follows that since $\osc_\beta f|_L = \osc_\alpha f|_L$ all
$\beta >\alpha$, $\sup_\beta \|\osc_\beta f|U\|_\infty = \|\osc_\alpha f|U
\|_\infty <\infty$, whence condition~4 of Theorem~5.1 holds, so
$f\in \DSC(K)$ by this result.\qed
\enddemo

Theorem 5.2 yields various solutions to the problem mentioned before its
statement.
The following corollary is a ``useable'' such solution.
Let us say that a series of $\real\cup \{\infty\}$-valued functions
$(g_j)$ on a Polish space $K$ converges {\it locally uniformly\/} if
for all $x\in K$, there is a neighborhood $U$ of $X$ and an $\nu$ so that
$g_j|U$ is real-valued for all $j\ge \nu$, and $\sum_{j\ge \nu} g_j$
converges uniformly on $U$.

\proclaim{Corollary 5.4}
Let $K$ be a Polish space, $(\varphi_j)$ a sequence in $\DSC(K)$, and $f$
a function on $K$ so that
\roster
\item"(i)" $\sum \varphi_j$ converges locally uniformly to $f$
\endroster
and
\roster
\item"(ii)" $\sum \osc_\alpha \varphi_j$ converges locally uniformly,
for every $\alpha <\omega_1$.
\endroster
Then $f\in \DSC(K)$.
\endproclaim

\demo{Proof}
Fix $x\in K$, $\alpha <\omega_1$, and choose $U$ a closed neighborhood of $X$
and a $\nu$ so that $\sum_{j=\nu+1}^\infty \varphi_j$ and
$\sum_{j=\nu+1}^\infty \osc_\alpha \varphi_j$ converge uniformly on $U$.
It follows that we may choose $\nu \le n_1<n_2<\cdots$ so that
$$\Big\| \sum_{j=n_i+1}^{n_{i+1}} \osc_\alpha \varphi_j |U\Big\|_\infty
< {1\over 2^i}\ \text{ for all } i=1,2,\ldots\ .
\tag 21$$
Now set $g= f-\sum_{j=1}^{n_1} \varphi_j$ and
$\psi_i = \sum_{j=n_i+1}^{n_{i+1}}\varphi_j $ for all $i$.
Then of course $\sum \psi_i$ converges uniformly on $U$ to $g$.
But moreover, for all $i$,
$$\osc_\alpha  \psi_i \le \sum_{j= n_i+1}^{n_{i+1}} \osc_\alpha \varphi_i
\ \text{ pointwise.}
\tag 22$$
Hence by (21) and (23),
$$\sum_{i=1}^\infty \|\osc_\alpha \psi_i|U\|_\infty <\infty\ .
\tag 23$$
It now follows by Theorem 5.2 that $g|U \in \DSC(U)$.
Indeed, for any non-empty closed subset $L$ of $U$, any $i$,
$\|\osc_\alpha \psi_i |L\cap U\|_\infty \le \|\osc_\alpha \psi |U\|_\infty$,
hence the hypotheses of Theorem~5.2 are fulfilled on $U$.
Then since $\DSC (U)$ is a linear space, also $f|U\in \DSC(U)$.
But this implies immediately that $f\in \DSC(K)$, since $f$ thus
locally belongs to $\DSC$.\qed
\enddemo

\enddocument